\documentclass[%
 aip,
 amsmath,amssymb,
 reprint,%
]{revtex4-1}

\usepackage{graphicx}
\usepackage{dcolumn}
\usepackage{bm}

\usepackage[utf8]{inputenc}
\usepackage[T1]{fontenc}
\usepackage{mathptmx}
\usepackage{etoolbox}
\usepackage{xcolor}

\usepackage[ruled,lined]{algorithm2e}

\def\HiLi{\leavevmode\rlap{\hbox to \hsize{\color{yellow!50}\leaders\hrule height .8\baselineskip depth .5ex\hfill}}}

\makeatletter
\def\@email#1#2{%
 \endgroup
 \patchcmd{\titleblock@produce}
  {\frontmatter@RRAPformat}
  {\frontmatter@RRAPformat{\produce@RRAP{*#1\href{mailto:#2}{#2}}}\frontmatter@RRAPformat}
  {}{}
}%
\makeatother
\begin{document}

\newcommand{\fk}[1]{{\color{blue}FK: #1}}

\preprint{AIP/123-QED}

\title[Truncated contagion maps]{Topological data analysis of truncated contagion maps}
\author{Florian Klimm}
\email{klimm@molgen.mpg.de}
\affiliation{ 
Department of Computational Molecular Biology, Max Planck Institute for Molecular Genetics, Ihnestra\ss{}e 63-73, D-14195, Berlin, Germany
}%
\affiliation{%
Department of Computer Science, Freie Universit\"at Berlin, Arnimallee 3, D-14195 Berlin, Germany
}%

\date{\today}

\begin{abstract}

 The investigation of dynamical processes on networks has been one focus for the study of contagion processes. It has been demonstrated that contagions can be used to obtain information about the embedding of nodes in a Euclidean space. Specifically, one can use the activation times of threshold contagions to construct contagion maps as a manifold-learning approach. One drawback of contagion maps is their high computational cost. Here, we demonstrate that a truncation of the threshold contagions may considerably speed up the construction of contagion maps. 
  Finally, we show that contagion maps may be used to find an insightful low-dimensional embedding for single-cell RNA-sequencing data in the form of cell-similarity networks and so reveal biological manifolds. Overall, our work makes the use of contagion maps as manifold-learning approaches on empirical network data more viable.

\end{abstract}

\maketitle

\begin{quotation}

It is known that the analysis of spreading processes on networks may reveal their hidden geometric structures.  These techniques, called contagion maps, are computationally expensive, which raises the question of whether they can be methodologically improved. Here, we demonstrate that a truncation (i.e., early stoppage) of the spreading processes leads to a substantial speed-up in the computation of contagion maps. For synthetic networks, we find that a carefully chosen truncation may also improve the recovery of hidden geometric structures. We quantify this improvement by comparing the topological properties of the original network with the constructed contagion maps by computing their persistent homology.  Lastly, we explore the embedding of single-cell transcriptomics data and show that contagion maps can help to distinguish different cell types. 
\end{quotation}

\section{Introduction}

The study of spreading processes has a rich tradition in epidemiology~\cite{kermack1991contributions,keeling2011modeling,schlosser2020covid}, the social sciences~\cite{valente1996social,christakis2013social,guilbeault2018complex,granovetter1978threshold,wiedermann2020network,oh2018complex}, and applied mathematics~\cite{borge2013cascading,porter2014dynamical,gleeson2008cascades}. Due to improved availability of interaction data, the description of contagion processes on networks has been a focus in recent years~\cite{watts1998collective}. It is known that spreading phenomena are influenced by the geometry of a networks's underlying embedding~\cite{brockmann2013hidden}, for example for the spread of pandemics~\cite{marvel2013small} and the adaption of technologies~\cite{ammerman1971measuring}. By constructing \emph{diffusion maps}, the interplay between dynamical processes and underlying geometry can be used to detect manifold structure in high-dimensional data~\cite{rohrdanz2011determination,ferguson2011nonlinear,coifman2008diffusion,kim1999nonlinear}. More recently, \emph{contagion maps} have been proposed as a technique to detect manifold structure in networks~\cite{taylor2015topological}. Contagion maps may reliable detect manifold structure in noisy data, even if standard approaches, such as {\sc isomap} fail~\cite{mahler2020contagion}. Generalisation of these approaches, for example to contagions on simplicial complexes~\cite{kilic2022simplicial} and Kleinberg-like networks~\cite{mahler2021analysis}, is an active field of research. The construction of contagion maps, however, is computationally expensive, which makes the exploration of strategies for their speed-up a pressing issue. In the context of random walks on networks, it has been established that stopping a dynamical process and restarting it may improve the performance of machine-learning algorithms~\cite{tong2006fast, kim2008generative,shun2016parallel}, raising the questions whether similar approaches might be able to enhance contagion maps.

\begin{figure*}[t]{
      \centering 
    \includegraphics[width=0.99\textwidth]{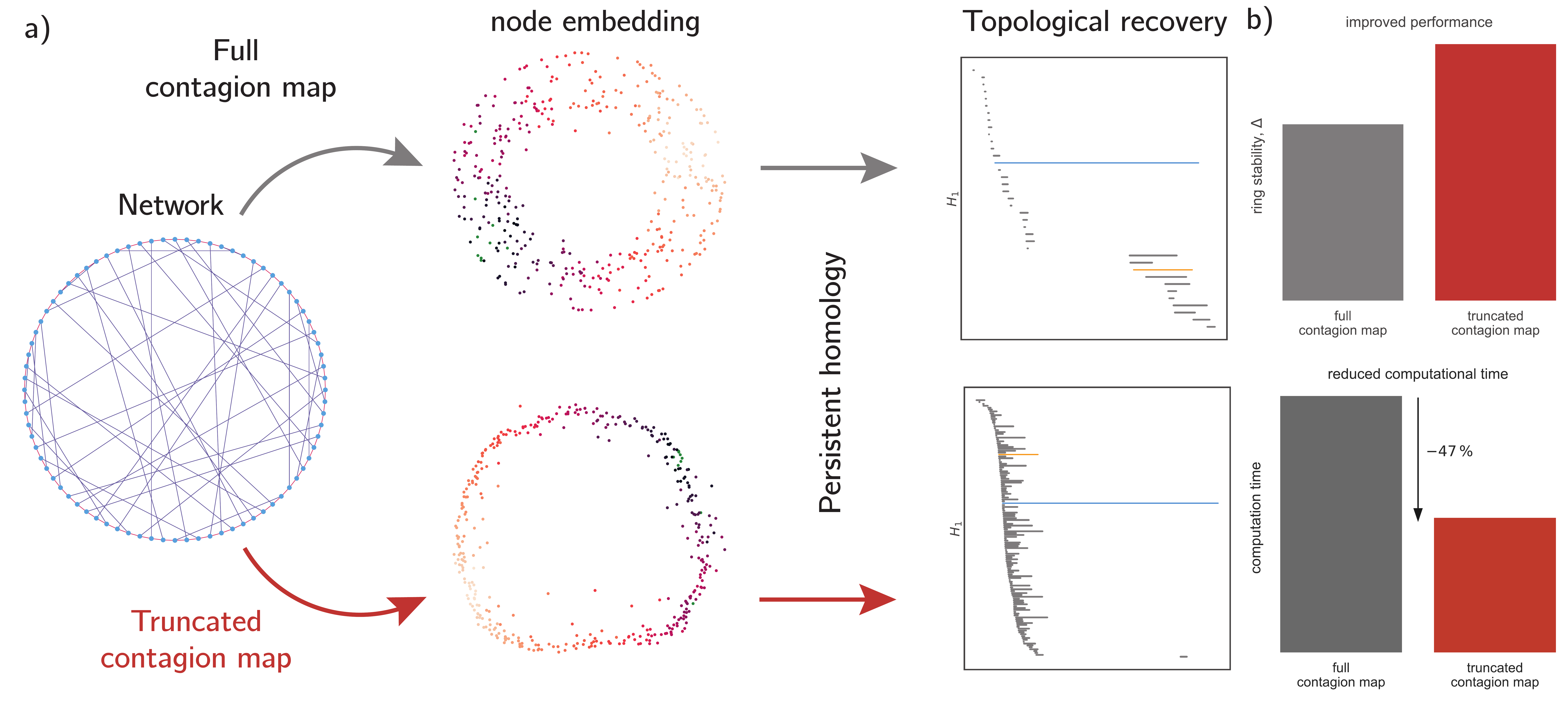}  
}
  \caption{\textbf{Truncated contagion maps can improve recovery of topological features of a networks' underlying manifold while having reduced computational cost.} (a) We show a schematic representation of the analysis pipeline in this manuscript: With contagion maps we can compute a node embedding for a network. With persistent homology, we compute barcode diagrams to quantify the recovery of topological features in the learned embeddings. This approach is used for the full contagion map and the newly introduced truncated contagion map. (b) We compare the topological recovery of full contagion maps with that of truncated contagion maps. For appropriate parameter choices, truncated contagion maps show higher performance than full contagion maps while taking considerably less time to compute (numbers for a noisy ring lattice with $N=1000$ nodes).
  }
  \label{fig:schema}
\end{figure*}

In this work, we introduce a generalisation of contagion maps that we call \emph{truncated contagion maps} (see Fig.~\ref{fig:schema}). The core idea is that an early stoppage of the contagion processes that are used to detect the manifold may drastically reduce the computational cost of the embedding algorithm. In particular, our numerical observations indicate a quadratic scaling in problem size for full contagion maps and a sub-quadratic scaling for truncated contagion maps. Furthermore, as we will show, the activations in the early stages of a contagion often follow the underlying manifold of the networks more clearly than the later stages. Omitting these late-stage dynamics results in an improved manifold detection, which we verify with techniques from \emph{topological data analysis}.

Topological data analysis is a field that aims to develop tools that allow the extraction of qualitative features in data that are hidden to standard approaches. A wide range of tools and methods for the topological analysis of data have been developed~\cite{otter2017roadmap,bianconi2015interdisciplinary,bick2021higher} and applied to problems in biology\cite{stolz2017persistent,giusti2016two}, social sciences~\cite{stolz2016topological,feng2021persistent,feng2020spatial}, and many others~\cite{stolz2020geometric,salnikov2018co,sizemore2018knowledge}. \emph{Persistent homology}~\cite{edelsbrunner2008persistent,edelsbrunner2013persistent} in particular is a topological data analysis method that received a lot of interest from practitioners because it enables the extraction of topological features in high-dimensional point-cloud data. In particular, topological data analysis is considering higher-order interactions between data points by constructing simplicial complexes.

Advances in \emph{single-cell transcriptomics} enable the measurement and analysis of gene expression at a single-cell resolution~\cite{aldridge2020single}. One challenge in the analysis of single-cell transcriptomics data is that experimental advances lead to an explosion in data set sizes~\cite{lahnemann2020eleven}. As complex organisms have thousands to ten thousands of genes, their gene-expression space is high-dimensional, making single-cell transcriptomics data notoriously hard to visualise. Yet, the vast majority of this space is empty and cells follow few low-dimensional manifolds that represent, for example, cell types or differentiation trajectories. The presence of low-dimensional structures embedded in high-dimensional gene-expression spaces makes single-cell transcriptomics a promising field for the application of manifold-learning approaches~\cite{moon2018manifold}. In this work, we apply truncated contagion maps to single-cell transciptomics data and show that the constructed embeddings are fruitful in representing the complex biological structure, which represents cell types and developmental trajectories.

The remainder of this article is organised as follows. In Section~\ref{sec:methods}, we summarise the existing methods and introduce truncated contagion maps. In Section~\ref{sec:bifurcation}, we investigate the temporal development of the contagion processes. In Section~\ref{sec:topo}, we demonstrate that truncated contagion maps may improve the recovery of topological features in contagion maps. In Section~\ref{sec:complexity}, we discuss the improved computational complexity of truncated contagion maps versus full contagion maps. In Section~\ref{sec:singleCell}, we apply truncated contagion maps to single-cell transcriptomics data of developing mouse oocytes. In Section~\ref{sec:conclusion}, we discuss our findings.

\section{Methods} \label{sec:methods}
\subsection{Noisy geometric networks}

Noisy geometric networks\cite{taylor2015topological} are geometric networks~\cite{barthelemy2011spatial} to which `noisy' non-geometric edges have been added. They are defined as follows. Consider a manifold $\mathcal{M}$ that is embedded in an ambient space $\mathcal{A}$ such that $\mathcal{M} \subset \mathcal{A}$. First, we construct nodes $V$ that are embedded in the manifold $\mathcal{M}$ such that $\mathbf{w}^{(i)} \in \mathcal{M}$ for nodes $i \in V$. Second, we construct edges which are of one of two types, either \emph{geometric} or \emph{non-geometric}. The geometric edges $E^{(G)}$ follow the manifold such that $i,j \in E^{(G)}$ if the distance along the manifold between nodes $i$ and $j$ is below some distance threshold. The non-geometric edges $E^{(NG)}$ are constructed in a random process that does not consider the manifold. We obtain the noisy geometric network as $G=(V,E^{(G)} \cup E^{(NG)})$. Noisy geometric networks are a type of small-world network~\cite{watts1998collective} which consist of long-range and short-range connections. The interplay of these two types of interactions has been discussed in the context of many types of networks, for example, social networks~\cite{granovetter1973strength} and traffic networks~\cite{brockmann2013hidden}.

\begin{figure}[b]{
      \centering 
  \includegraphics[width=0.48\textwidth]{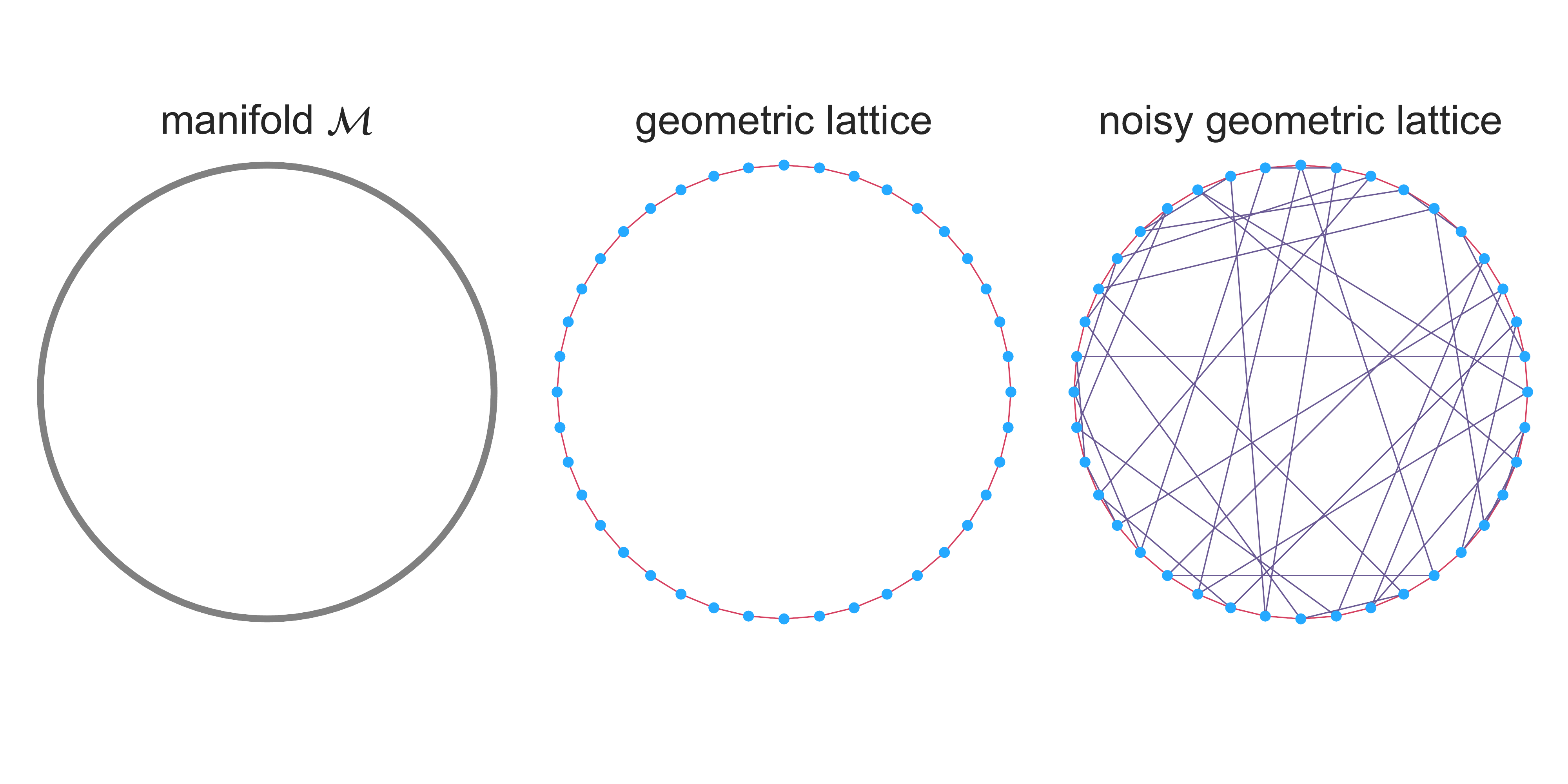}
  }
  \caption{\textbf{Noisy geometric networks contain geometric edges that follow the underlying manifold and non-geometric edges that ignore the underlying manifold.}  Here, we show noisy geometric lattices of size $N=50$, whose underlying manifold $\mathcal{M}$ is the unit circle $\{ (a,b) | a^2 +b^2=1\}$ with ambient space $\mathcal{A} \subset \mathbb{R}^2$. The geometric lattice contains exclusively geometric edges and has $d^{(G)}=2$. After adding $d^{(NG)}=2$ non-geometric edges per node, we obtain a noisy geometric lattice with $\alpha=1$.}
  \label{fig:noisyGeometric}
\end{figure}

In this manuscript, we focus on noisy ring lattices $NRL(N,d^{(G)},d^{(NG)})$ as prototypical examples of noisy geometric networks. Their parameters are the number $N$ of nodes, the geometric degree $d^{(G)}$, and the non-geometric degree $d^{(NG)}$. The NRL networks are defined as follows. First, we place $N$ nodes uniformly spaced along the unit circle $\mathcal{M} = \{ (a,b) | a^2 +b^2=1\} \subset \mathbb{R}^2$ such that $\mathbf{w}^{(i)}=(\cos (2\pi i/N),\sin (2\pi i/N ))$ is the location of node $i$. Second, we construct geometric edges such that each node is connected to its $d^{(G)}$ nearest-neighbour nodes. Third, we construct $N\times d^{(NG)}$ non-geometric edges uniformly at random such that each node is has exactly $d^{(NG)}$ non-geometric edges. In practice, the non-geometric edges are constructed with \emph{stub-matching}: We create a list of stubs that contains each node exactly $d^{(NG)}$ times. Then we construct edges iteratively by drawing without replacement uniformly at random from this list while rejecting self-edges or parallel edges. In rare cases this procedure might fail to terminate (e.g., if only stubs of the same node are remaining) but this problem can be overcome by restarting the procedure. The ratio $\alpha=d^{(NG)} / d^{(G)}$ of non-geometric to geometric edges quantifies the `noisiness' of a noisy ring lattice with $\alpha =0$ indicating that there are exclusively geometric edges, and $\alpha =1$ indicating that there are as many non-geometric edges as geometric edges. In Fig.~\ref{fig:noisyGeometric}, we show two example noisy geometric lattices with $\alpha=0$ and $\alpha=1$, respectively. 

\subsection{Watts threshold model}

The Watts threshold model (WTM) is a well-established deterministic binary-state discrete-time model of social contagions on networks~\cite{watts2011simple} and is a modification of Mark Granovetter's threshold model~\cite{granovetter1978threshold}. The WTM is defined as follows. Each node is in one of two states, either active or inactive. At time $t=0$ all nodes but a set of \emph{seed nodes} $S$ is inactive. All nodes have the same activation threshold $T\in [0,1]$ (homogenous threshold), which determines how easily a node adapts their neighbours' activation. An active node stays active but an inactive node is activated if more than a threshold fraction $T$ of its $d$ neighbours are active.  We update all node states synchronously at each time point and stop the contagion once a steady state is reached (for a discussion of synchronous vs asynchronous update procedures, see \cite{porter2016dynamical}). 

\subsection{Contagion maps}

As introduced by \citeauthor{taylor2015topological}, we use the WTM to construct deterministic embeddings for networks. These embeddings are based on the nodes' activation times under different initial conditions. Specifically, we define a contagion map as $V \to \{\mathbf{x}^{(i)}\}_{i \in V}$, where $\mathbf{x}^{(i)} = [x_1^{(i)},x_2^{(i)},\dots,x_j^{(i)}]$ with $x_j^{(i)}$ indicating the activation time of node $i$ for contagion $j$. For the initial conditions, we use \emph{cluster seeding} such that the neighbourhood $\mathcal{N}(j)$ of node $j$ is active at $T=0$. After a steady state of the WTM is reached some nodes might still be inactive. The infinite activation time of these inactive nodes is set to $2N \ll \infty$ because $2N$ is much larger than the largest possible activation time of $N-1$ (An activation time of $N-1$ would be reached by a single node if in each time step exactly one additional node is activated, for example, along a line graph.)

The obtained matrix $\mathbf{X}$ of activation times is not necessarily symmetric. To study an embedding that defines a semimetric, we investigate the symmetric contagion map $\mathbf{X}_{\text{symmetric}} = \mathbf{X} + \mathbf{X}^\top$. See Algorithms~\ref{alg:contagionMaps} and \ref{alg:WTM} for pseudocode for the contagion maps and contagion dynamics, respectively.

For clarity, we will refer to the contagion maps as defined by \citeauthor{taylor2015topological} as 'full contagion maps' in contrast to the `truncated contagion maps' that we will introduce now.

\begin{algorithm}
\caption{Computing (truncated) contagion maps}\label{alg:contagionMaps}
         \SetKwInOut{Input}{input}
        \SetKwInOut{Output}{output}
        \Input{Network $G=(V,E)$ with $N=|V|$ nodes\\
        Threshold $T$\\
        Number of steps $s$ \tcp{Ignored for full contagion maps}}
        \Output{Node embedding $\mathbf{X}_{\text{symmetric}} \in \mathbb{N}^{N\times N}$}
    \SetKwBlock{Beginn}{beginn}{ende}
    \Begin{
        \For{each node $j \in V = \{1, \dots, N\}$}{
                Initialise cluster seeding by $S$ infecting $j$ and its neighbours\;
				$\mathbf{X[:,j]} \leftarrow $ Compute (truncated) contagion dynamics\;
        }
        $\mathbf{X}_{\text{symmetric}} = \mathbf{X} + \mathbf{X}^\top$\;
        \Return $\mathbf{X}_{\text{symmetric}}$
    }
\end{algorithm}

\begin{algorithm}
\caption{Compute (truncated) contagion dynamics} \label{alg:WTM}
         \SetKwInOut{Input}{input}
        \SetKwInOut{Output}{output}
        \Input{Network $G=(V,E)$ with $N=|V|$ nodes\\seed nodes $S\subset V$\\ threshold $T$\\max number of steps $s$}
        \Output{Activation time $\vec{x} \in \mathbb{N}^{N}$}
    \SetKwBlock{Beginn}{beginn}{ende}
    \Begin{
    	Initialise activation times $\vec{x} = \vec{\text{NaN}}$\\
    	Initialise step number $t=0$\\
    	\While{$\vec{x} \neq \vec{x}$ }{
    	\For{each node $i \in V = \{1, \dots, N\}$}{Check if node $i$ activated\\
    	Update $\vec{x}$
        }
        \tcp{Break condition only checked for truncated contagion maps}
        \If{$t==s$}{
        Set activation time to $2\times s$ for inactive nodes\\
        break while loop} 
        }
        
        \Return $\vec{x}$
    }
\end{algorithm}

\subsection{Truncated contagion maps}

In this work, we introduce \emph{truncated contagion maps} as a generalisation of full contagion contagion maps. Intuitively, they represent full contagion maps in which contagions are not run until a steady state is reached. Rather, we stop the contagion after $s \in \mathbb{N}_{>0}$ steps. Specifically, we define a truncated contagion map as $V \to \{\mathbf{x}^{(i)}\}_{i \in V}$, where $\mathbf{x}^{(i)} = [x_1^{(i)},x_2^{(i)},\dots,x_j^{(i)}]$ with

\begin{align}
	x_j^{(i)} = \begin{cases}
a_j^{(i)} & \text{for}\ a_j^{(i)}\leq s\\
2s & \text{else}\,.
\end{cases}
\end{align}

where $a_j^{(i)}$ indicates the activation time of node $i$ for contagion $j$. For the initial conditions, we use \emph{cluster seeding} such that the neighbourhood $\mathcal{N}(j)$ of node $j$ is active at $T=0$. We set the activation time of unactivated nodes to $2s$, such that full contagion maps are a special case of truncated contagion maps for $s=N$, but other choices for the handling of infinite activation times are possible. In practice, however, we anticipate an approximately similar behaviour of truncated and full contagion maps for $s\ll N$ because most contagions reach a steady state after a small number $s$ of steps (in comparison to the network size $N$). As for full contagion maps, we study exclusively the symmetric contagion maps $\mathbf{X}_{\text{symmetric}} = \mathbf{X} + \mathbf{X}^\top$ such that we obtain a semimetric.

\subsection{Topological data analysis}

Persistent homology is a widely adapted method from topological data analysis. It aims to identify 
stable (i.e., persistent) topological features of point cloud data across a range of resolutions. To obtain these topological features, we must construct simplicial complexes from the point cloud data. Specifically, we employ an Vietoris--Rips simplicial filtration, which is a sequence of simplicial complexes that constructs topological structure based on the data. Given a data set $X$, a metric $d: Y \times Z \to \mathbb{R}$, and a scale parameter $\alpha \in \mathbb{R}$, we can construct a Vietoris--Rips complex as
\begin{align}
\{{VR}(X, d, \alpha) 	:= \{ \sigma \subset X: d(p,q) \leq \alpha\ \forall\ p,q \in \sigma \}\,. 
\end{align}
Computing the Vietoris--Rips complex for a selection of scale parameters $\alpha_1 \leq \alpha_2 \leq \dots \leq  \alpha_a$, then yields a Vietoris--Rips simplicial filtration 

\begin{align}
\emptyset \subseteq \text{VR}(X,d,\alpha_1) \subseteq \text{VR}(X,d,\alpha_2) \subseteq \dots \subseteq \text{VR}(X,d,\alpha_a) \,, 
\end{align}
which is a sequence of embedded simplicial complexes. This filtration can be though of as a dynamical process on a simplicial complex in which faces are added stepwise. The most common way to visualise the persistence of topological features across a filtration is a \emph{barcode}, which indicates in which filtration step a topological feature is born and in which step it ceases to exist.  We use the {\sc Ripser} implementation~\cite{tralie2018ripser} for the computation of Vietoris--Rips persistence barcodes and use the Euclidean norm. To summarise the 1D persistence barcodes (i.e., the presence of loops in the data), we compute \emph{ring stability} as $\Delta = l_1-l_2$, where $l_1$ and $l_2$ are the persistence of the longest and second longest 1ifetime of loops in the data. This procedure has previously been used to identify whether contagions are predominantly spreading by wavefront propagation versus new cluster appearance~\cite{taylor2015topological}. In Fig.~\ref{fig:schema}, we show two example Vietoris--Rips persistence barcodes constructed from contagion maps.

\section{Temporal development of contagions reveals suitable truncation values}
\label{sec:bifurcation}

\begin{figure}[t]{
      \centering 
  \includegraphics[width=0.45\textwidth]{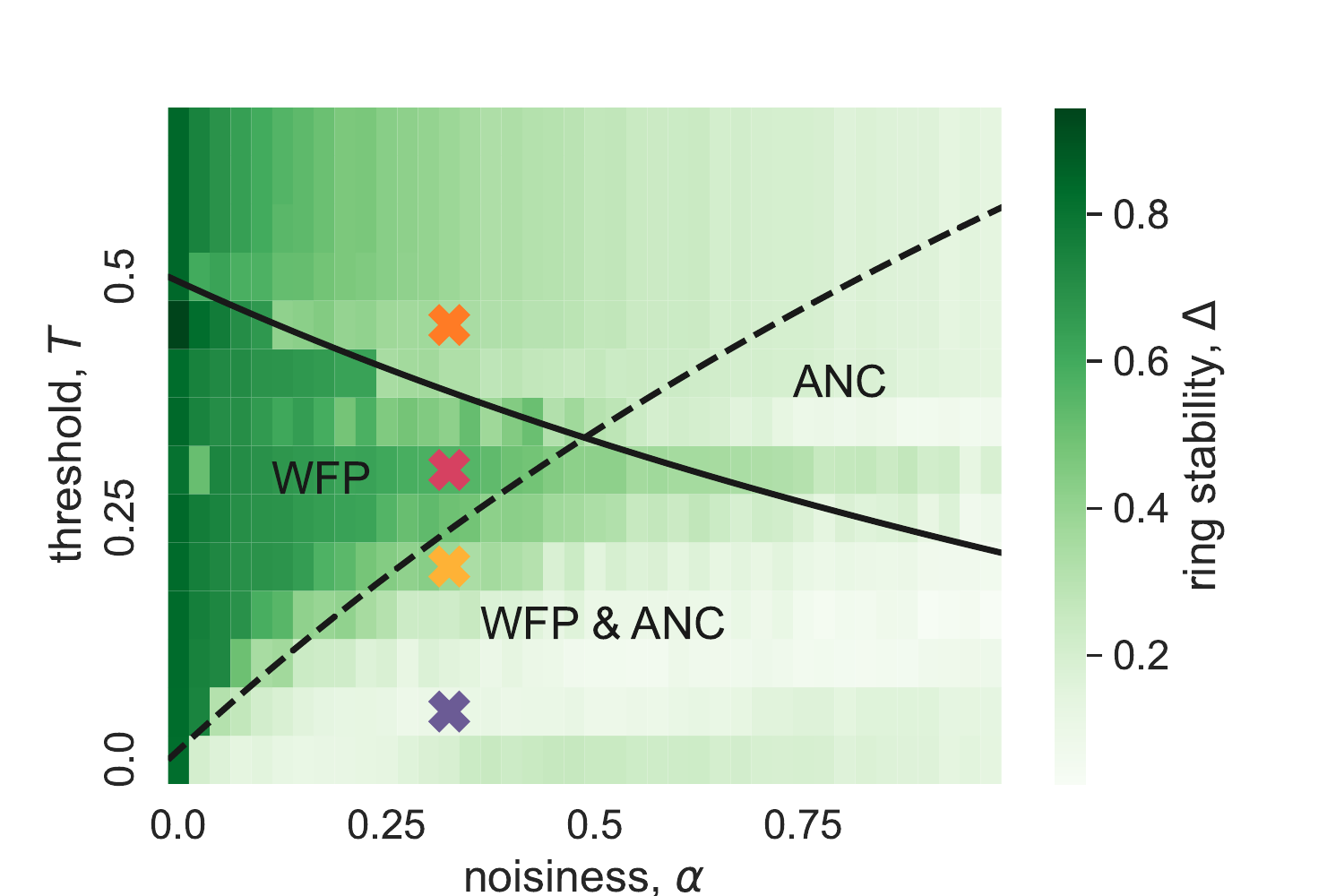}
  }
  \caption{\textbf{Bifurcation diagram of WTM on noisy ring lattices matches well the observed ring stability of truncated contagion maps.} We plot the critical thresholds of WFP (solid line, Eq.~\ref{eq:WFP}) and ANC (dashed line, Eq.~\ref{eq:ANC}) versus the noisiness $\alpha$. These curves divide the $(T,\alpha)$ parameter space into four qualitatively different contagion regimes: exclusively WFP, exclusively ANC, WFP and ANC, and neither.  
   We compare the theoretical bifurcation diagram with the empirically observed ring stability $\Delta$ in truncated contagion maps with $s=20$ of noisy ring lattices of size $N=400$ with geometric degree $d^{(G}=40$ and varying non-geometric degrees $d^{(NG} \in [0,40]$ such that their noisiness $\alpha \in [0,1]$. We find that only in the parameter regime of exclusive WFP the ring stability $\Delta$ is high, indicating a good recovery of the underlying manifold's topology. We highlight the four parameter choices that we show in Figs.~\ref{fig:temporalDevelopment} and~\ref{fig:embeddedMaps} with crosses.
  }
  \label{fig:matrixPlot}
\end{figure}

\begin{figure}[t]{
      \centering 
  \includegraphics[width=0.45\textwidth]{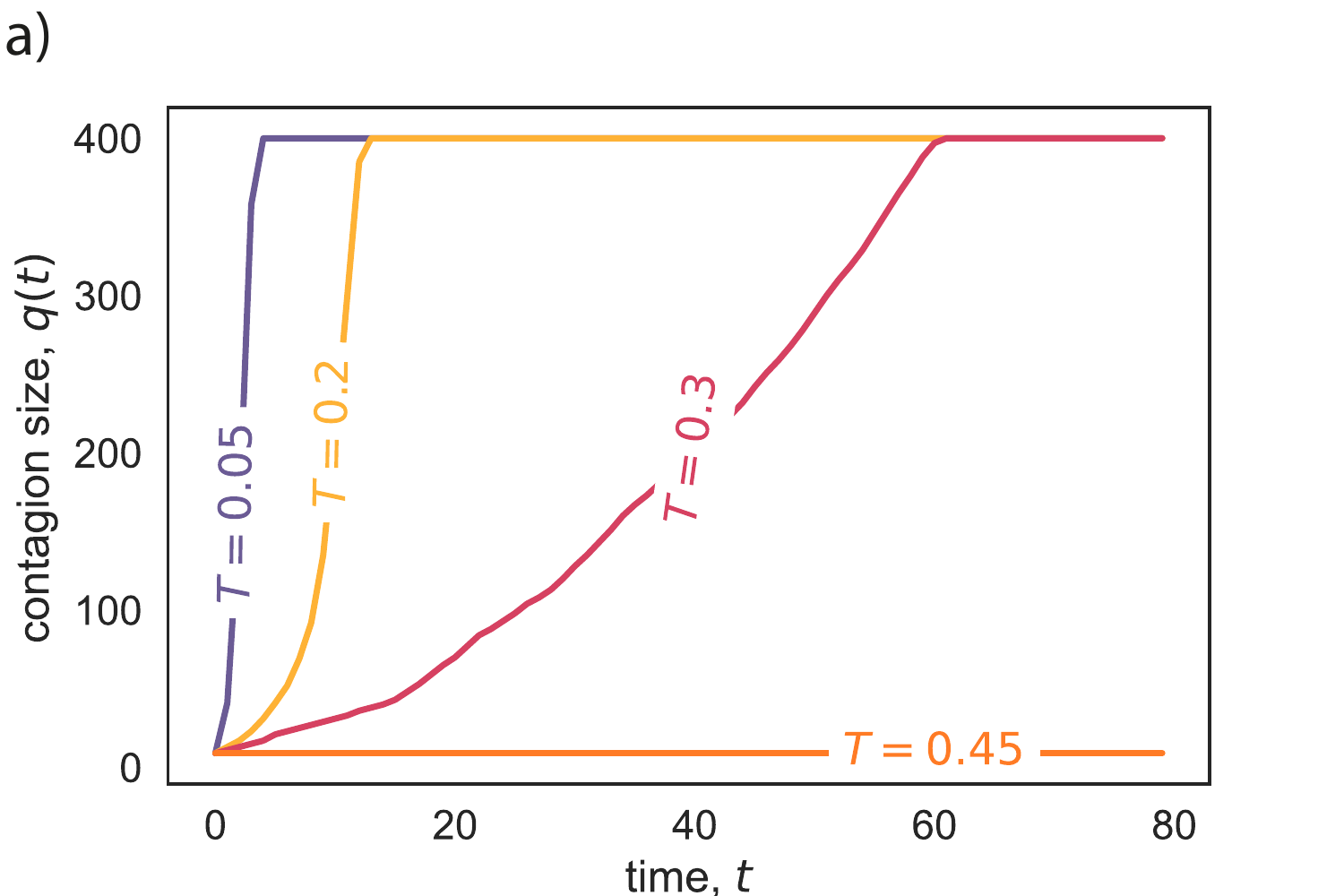}
  \includegraphics[width=0.45\textwidth]{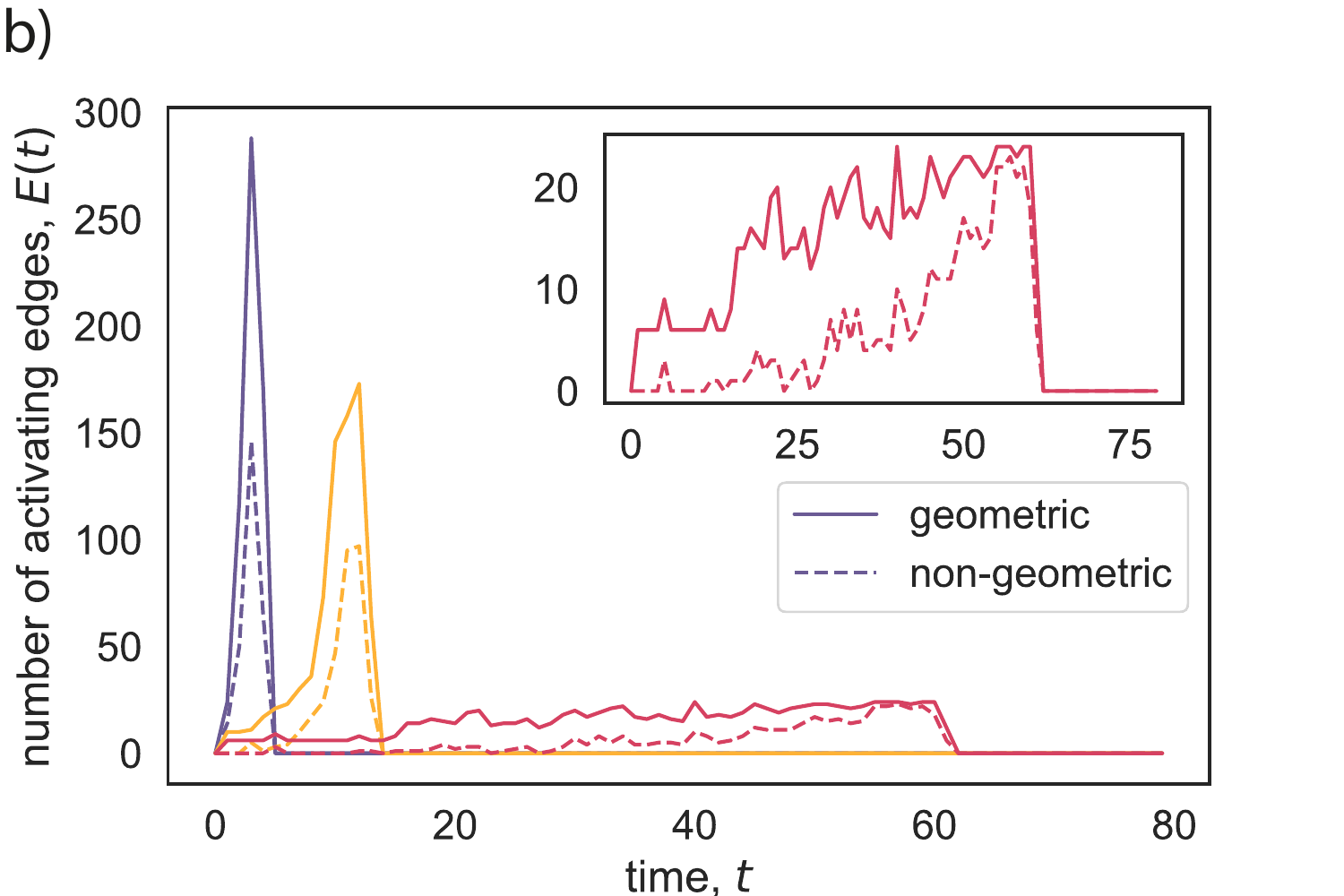}
  }
  \caption{\textbf{The temporal development of the contagion is strongly dependent on the threshold $T$.} (a) The size $q(t)$ of the contagion is growing over time for thresholds $T \in \{0.05,0.2,0.3\}$ and reaches $q(t)=1$, which indicates a global contagion (i.e., all nodes are active). The size $q(t)$ of the contagion remains constant for $T=0.45$, indicating that no cascade is triggered. (b) We show the number $E(t)$ of edges involved in the activation process and distinguish between geometric edges (solid lines) and non-geometric edges (dashed lines). The inlay shows a zoom onto the contagion with threshold $T=0.3$ in which we observe that for the early (i.e., $t<20$) contagions steps almost no non-geometric edges are involved. 
  }
  \label{fig:temporalDevelopment}
\end{figure}

In this section, we will describe the spread of WTM contagions on noisy ring lattices. In particular, we will demonstrate that in certain parameter regimes, the initial spread is mainly occurring along the manifold (called `wave front propagation', WFP) and at later stages of the contagion the spread is less restricted to the manifold and leads to new contagion clusters at distant parts of the manifold (`appearance of new clusters', ANC).

\citeauthor{taylor2015topological}~\cite{taylor2015topological} proved two critical thresholds for noisy ring lattices of noisiness $\alpha$: 

\begin{align}
	T^{\text{(WFP)}}&=1/(2+2\alpha)\,,\ \text{and} \label{eq:WFP}\\
	T^{\text{(ANC)}}&=\alpha/(1+\alpha)\,, \label{eq:ANC}
\end{align}
which indicate the thresholds above which no WFP and no ANC, respectively, occur. We show these critical thresholds, which intersect at $(\alpha,T)=(1/2,1/3)$, in Fig.~\ref{fig:matrixPlot} and observe four qualitatively different contagion regimes: exclusively WFP, exclusively ANC, WFP and ANC, and neither. As we will show, these regimes generalise from full contagion maps to truncated contagion maps. In the following, we investigate four parameter choices in more detail. Specifically, we summarise the Watts threshold model contagions for thresholds $T \in \{0.05,0.2,0.3,0.45\}$ on a noisy ring lattice with $N=400$ nodes and noisiness of $\alpha=d^{(NG)} / d^{(G)} = 2/6 = 1/3$ (see Fig.~\ref{fig:temporalDevelopment}).

In Fig.~\ref{fig:temporalDevelopment}a, we show the size $q(t)$ of the contagion over time (i.e., the number of nodes that are active). First, we note that for $T \in \{0.05,0.2,0.3\}$ a global cascade is triggered as $q(t) \to 1$, whereas for $T=0.45$ the contagion size remains constant. Second, a smaller threshold $T$ results in a quicker contagion spread which in turn yields a smaller time until a global cascade is reached. By comparing the chosen thresholds with the critical thresholds of WFP (Eq.~\ref{eq:WFP}) and ANC (Eq.~\ref{eq:ANC}) we observe that $T=0.05$ and $T=0.2$ fall in the regime in which WFP and ANC occur, $T=0.3$ falls in the regime in which exclusively WFP occurs, and $T=0.45$ results in no contagion (see coloured crossed in Fig.~\ref{fig:matrixPlot}).

To study whether the contagion spreads along the manifold or independent of it, we investigate along which type of edges the contagion spreads predominantly. In Fig.~\ref{fig:temporalDevelopment}b, we show the number $E(t)$ of activating edges (i.e., edges incident to a node that is active at $t$ and a node that becomes active at $t+1$). We distinguish between geometric edges (solid lines) and non-geometric edges (dashed lines). We observe that for thresholds $T \in \{0.05,0.2\}$ the contagion spreads quickly and to similar extent along geometric and non-geometric edges. For a threshold of $T=0.3$ in contrast, we find a slow contagion process. We also find that initially the contagion process is strongly dominated by contagion along geometric edges (see inlay). After time $t\approx 20$, however, the influence of non-geometric edges increases. Our observation indicates that for slow contagion processes, the contagion initially spreads along the manifold as WFP. Yet, a spreading along non-geometric edges becomes more important over time. This highlights that the derived critical thresholds (Eqs.~\ref{eq:WFP} and~\ref{eq:ANC}) are strictly valid only for the early contagion steps and become less appropriate the longer the contagion progresses. These observations indicate that a truncation choice that removes these late-stage contagion steps might improve the contagion map. In this case, for example, a truncation parameter $s\in [20,60]$ might be appropriate. As a simple, heuristic choice, we suggest $s^{\ast}=N/10=40$, which we will explore in more detail in Sec.~\ref{sec:topo}. As we will show in Sec.~\ref{sec:topo}, such truncations allow us to improve the recovery of topological features with truncated contagion maps in comparison with full contagion maps.

\section{Topological data analysis of truncated contagion maps}
\label{sec:topo}

\begin{figure*}[t]{
      \centering 
  \includegraphics[width=0.99\textwidth]{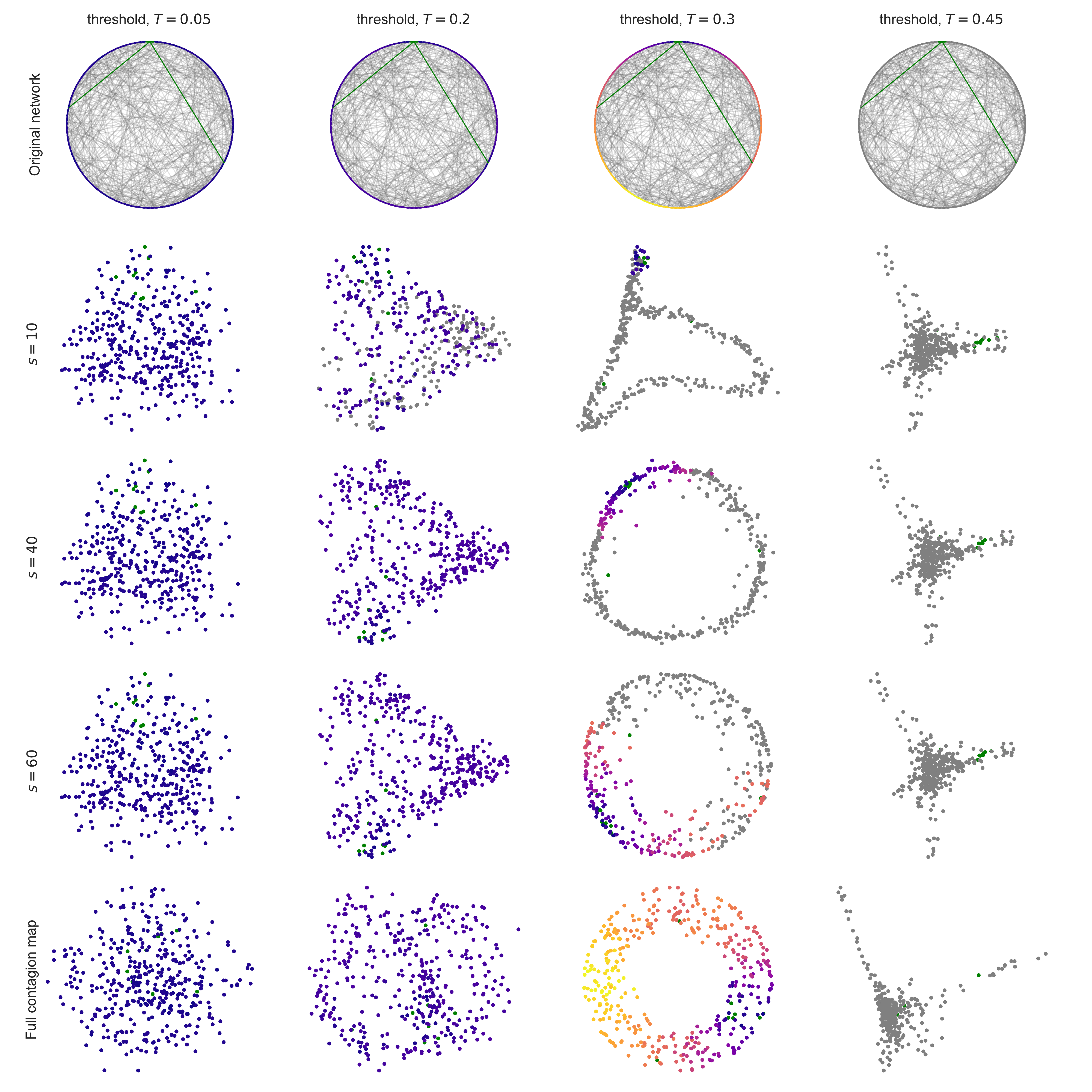}
  }
  \caption{\textbf{Truncated contagion maps applied to noisy ring lattices can recover the ring manifold for appropriate thresholds $T$}. We applied contagion maps to noisy ring lattices with $N=400$ nodes degrees $ = (d^{(G)},d^{(NG)}) = (6,2)$. In each column, we show the results for one of the four thresholds $T \in \{0.05,0.2,0.3,0.45\}$. (top row) The network with nodes embedded on the ring manifold. Each node's colour indicates the activation time from small (blue) to large (yellow) in one realisation of the WTM. Green nodes are seed nodes and grey nodes are inactive at the steady state. (other rows) We show two-dimensional projections of the $N$-dimensional contagion maps. A node's colour indicate its activation time as in the top row. The rows represent truncated contagion maps with $s=10$, $s=40$, and $s=60$, and a full contagion map, respectively. We observe that the threshold $T=0.3$ the contagion maps most closely resemble the circular embedding manifold. In this regime, intermediate step numbers $s=40$ leads to a best recovery of the circular manifold (up to rotation).
  }
  \label{fig:embeddedMaps}
\end{figure*}

\begin{figure*}[t]{
      \centering 
    \includegraphics[width=0.32\textwidth]{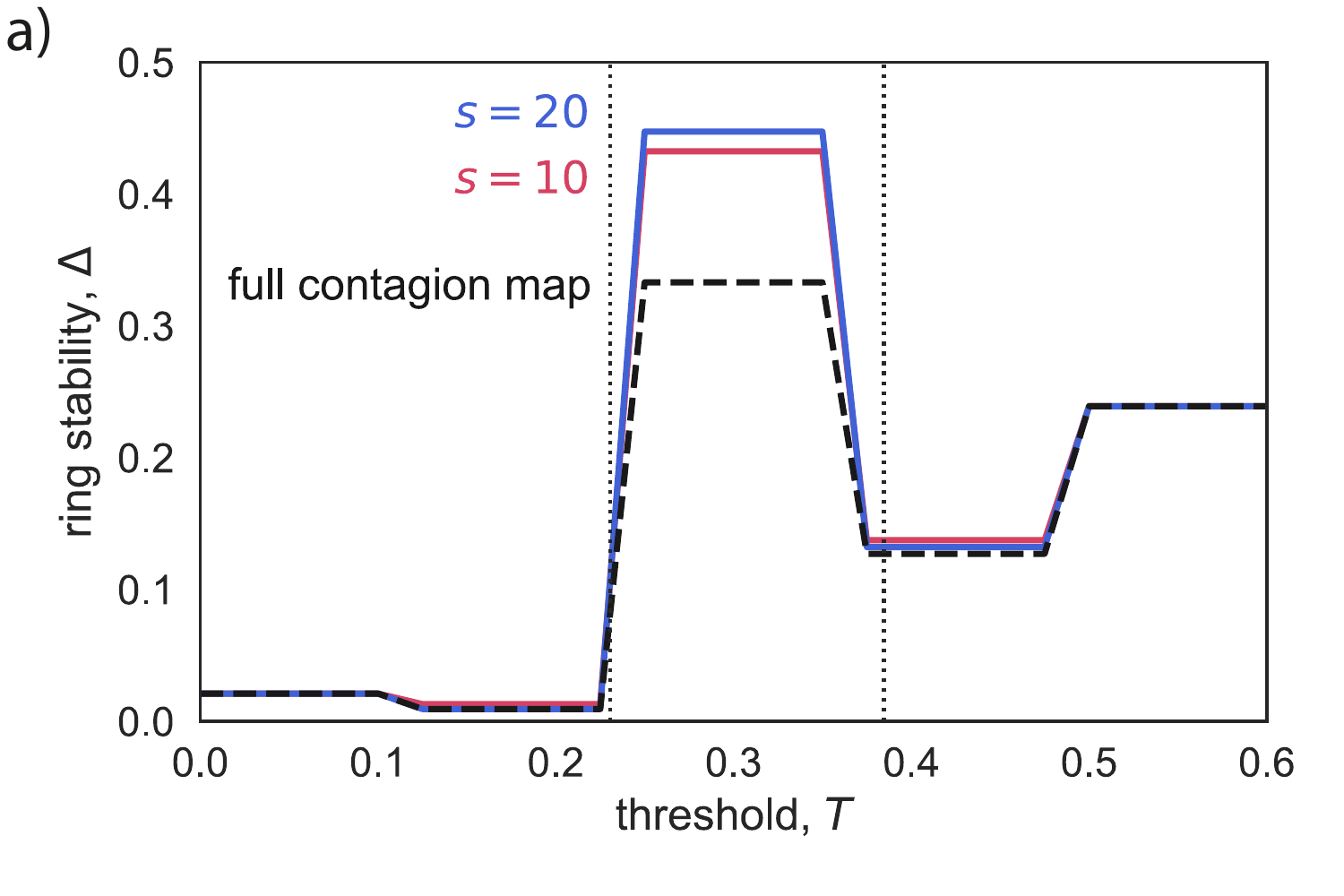}  
    \includegraphics[width=0.32\textwidth]{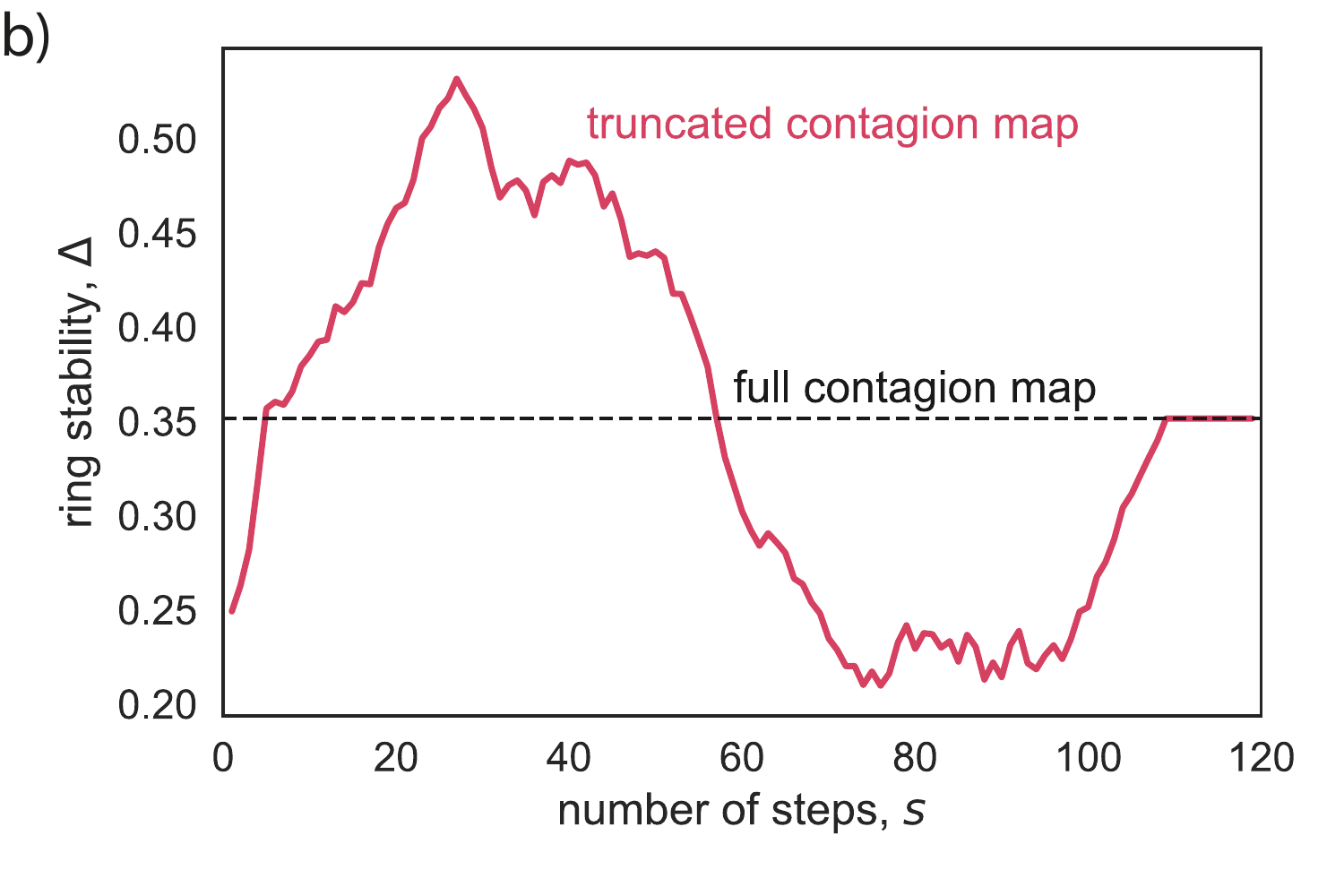}
    \includegraphics[width=0.32\textwidth]{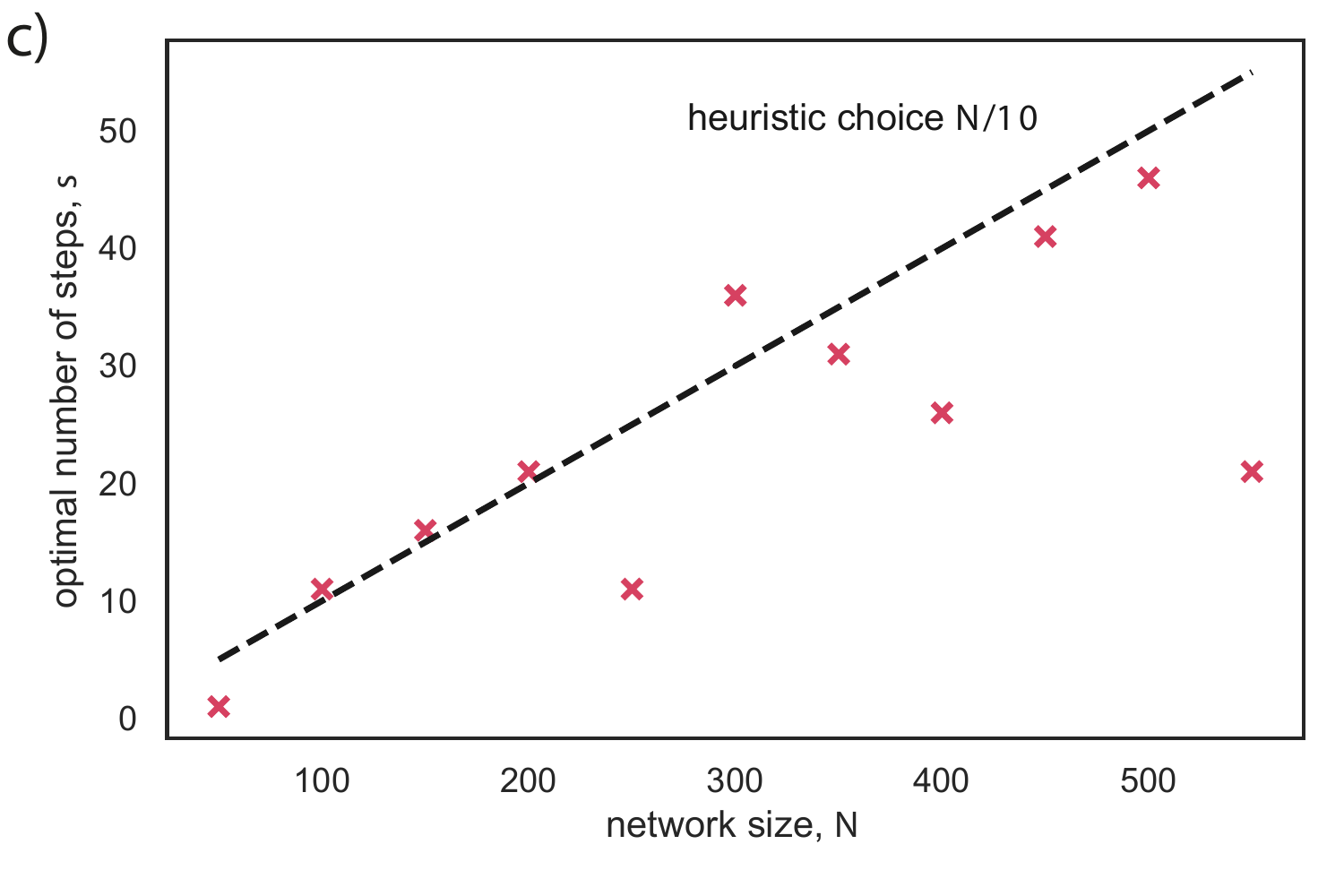}
}
  \caption{\textbf{Truncated contagion maps enable an improved recovery of topological features, in comparison with a full contagion map.} We quantify the recovery of the topological features by computing the ring stability $\Delta$ of contagion maps on noisy ring lattices with $N=400$. (a) For all contagion maps, the ring stability is largest in the threshold regime $\frac{\alpha}{1+\alpha} <T <\frac{1}{2+2\alpha}$ in which WFP dominates the spreading process. The truncated contagion maps with $s=10$ (red line) and $s=20$ (blue line) have a larger ring stability $\Delta$ than the full truncation map (dashed black line). (b) The ring stability $\Delta$ varies with the number $s$ of included contagion steps. For small number $s$ of contagion steps, the truncated contagion map results in an improved ring stability $\Delta$ in comparison for the full contagion map (dashed black line). For $s>60$, when ANC occurs, the ring stability $\Delta$ diminishes. As expected the performance of the full contagion map matches the truncated contagion map for large $s$. Results for the threshold $T=0.3$ (c) We investigate the truncation parameter $s^{\ast}$ that yields the best recovery of the ring manifold for noisy ring lattice networks of sizes $N \in [50,550]$. We find that the optimal truncation increases approximately linearly with network size $N$, which suggests a simple heuristic choice of $s^{\ast}=N/10$ (dashed line).
  }
  \label{fig:topologicalFeatures}
\end{figure*}

In this section, we analyse the topological properties of contagion maps on noisy ring lattices. This allows us to identify parameter regimes in which the contagion processes follow the underlying manifold. In particular, we compare the behaviour of truncated contagion maps with that of full contagion maps.

In Fig.~\ref{fig:topologicalFeatures}, we study contagion maps for noisy ring lattices with $N=400$ nodes and noisiness of $\alpha=d^{(NG)} / d^{(G)} = 2/6 = 1/3$. Each column represents a distinct threshold  $T \in \{0.05,0.2,0.3,0.45\}$, which correspond to the example thresholds in Fig.~\ref{fig:temporalDevelopment} and distinct regimes in the bifucation analysis of the critical thresholds critical thresholds (Eqs.~\ref{eq:WFP} and ~\ref{eq:ANC}). In the top row, we show the noisy ring lattice with nodes placed at their locations $\mathbf{w}^{(i)}=(\cos (2\pi i/N),\sin (2\pi i/N )$. The colour of each node indicates the node's activation time during one realisation of the WTM in which the green nodes are the seed nodes. For all other nodes, the colour map ranges from blue (small activation time) to yellow (large activation time). Grey nodes never adopt the contagion. As expected, we observe that for two smallest thresholds $T \in \{0.05,0.2\}$ the contagion quickly spreads to the whole network and that for the largest threshold $T=0.45$ no contagion beyond the seed nodes occurs. For intermediary threshold $T=0.3$, however, we observe that the contagion spreads approximately along the circular manifold.

In the remaining rows of Fig.~\ref{fig:topologicalFeatures}, we show the point clouds $\mathbf{x}^{(i)}$ obtained from contagion maps. Specifically, rows two, three, and four show truncated contagion maps for $s=10$ steps, $s=40$ steps, and $s=60$ steps and the last row shows the full contagion map, which is equivalent to $s=N$. To visualise the $N$-dimensional point clouds, we use a two-dimensional projection through principal component analysis (PCA). Similar to the networks in the top row, we colour each point with the activation time under one realisation of the WTM. Yet, as we truncate the contagion after $s$ steps a node might be grey (i.e., never activated) either because a steady state is reached or because the contagion stopped earlier. Therefore, the number of grey (i.e., never activated) nodes varies with the number $s$ of considered contagion steps although the seed nodes were identical.


In general, we find that the contagion maps $\mathbf{x}^{(i)}$ for threshold $T=0.3$, which is in a regime in which there is exclusively WFP predicted, most closely resemble the two-dimensional ring topology of the original network (up to rotation). Furthermore, in this regime, the truncation parameter $s$ has a strong influence on the obtained embedding. For a small number $s=10$ of contagion steps, the ring structure is perturbed and appears folded. Increasing the number of contagion steps to $s=20$ makes the ring structure more clearly visible. Increasing the number $s$ of contagion steps further, however, results in a `broadening' of the circular point clouds. Our observation suggest that for the threshold $T=0.3$ a truncated contagion map might better recover the underlying circular manifold than the full contagion map.

After the bifurcation analysis in Sec.~\ref{sec:bifurcation} and the qualitative investigation in Fig.~\ref{fig:embeddedMaps}, we now quantify the recovery of the manifolds topological features with persistent homology. In particular, we know that the underlying manifold $\mathcal{M}$ of a noisy ring lattice is the unit circle. By comparing how closely the persistent homology of the contagion maps resemble the topological properties of the unit circle, we can quantify the recovery of its topological properties. A unit circle has one circular hole and thus its first Betti number is $b_1=1$.  Therefore, we compute the ring stability $\Delta$ as a quantification of how `ring-like' a contagion map is. In Fig.~\ref{fig:topologicalFeatures}a, we show the ring stability $\Delta$ of contagion maps on noisy ring lattices of size $N=400$ and noisiness $\alpha=1/3$. We vary the threshold $T\in [0,0.6]$ and highlight the critical thresholds $T^{\text{(ANC)}} = 0.25$ and $T^{\text{(WFP)}} = 0.375$ (Eqs.~\ref{eq:WFP} and ~\ref{eq:ANC}) as vertical dotted lines. All contagion maps have their largest ring stability $\Delta$ in the regime between these thresholds (i.e., when only WFP is possible). Furthermore, the truncated contagion maps with $s=10$ and $s=20$ have a higher ring stability $\Delta$ than the full contagion map, which matches our observations in Fig.~\ref{fig:embeddedMaps} that an appropriate truncation of the contagion map may improve the recovery of the ring-like topology. In other threshold regimes, the ring stability $\Delta$ is almost indistinguishable between truncated and full contagion maps.

In Fig.~\ref{fig:topologicalFeatures}b, we investigate the influence of the number $s$ of contagion steps on the ring stability $\Delta$ of contagion maps with a threshold of $T=0.3$. We find that for a large number $s\to N$ of steps, the ring stability is the same as for the full contagion map (dashed line), which matches our theoretical expectation (see `Truncated contagion maps' in the Methods section). For a small number $s<60$ of steps, the truncated contagion map has a larger ring stability $\Delta$ than the full contagion map, whereas for $s>60$ steps, the full contagion map outperforms the truncated contagion map. We explain this behaviour as a trade-off between data quality and data quantity: The early-stage contagions more closely follow the underlying manifold, yielding better embeddings for the truncated contagion maps. Yet, full contagions have more data as less entries in the contagion map have infinite activation times.

The dependence of the ring stability $\Delta$ on the truncation parameter $s$ raises the question whether  the optimal choice varies with the network size $N$. In Sec.~\ref{sec:bifurcation}, we suggested $s^{\ast} = N/10$ as a simple, heuristic choice because in this regime, the contagion is predominantly spreading along then manifold. In Fig.~\ref{fig:topologicalFeatures}c, we compare this heuristic with the actual optimal $s^{\ast}$ for noisy ring-lattice networks of sizes $N \in [50,550]$. We find that the optimal truncation parameter $s^{\ast}$ is increasing with the network size $N$ and is in a decent agreement with the suggested heuristic. 

Thus far, we fixed the noisiness $\alpha$ of the networks. In Fig.~\ref{fig:matrixPlot}, we compare the analytically predicted bifurcation diagram with the ring stability $\Delta$ of truncated contagion maps with $s=20$ steps for various thresholds $T \in [0,0.6]$ for noisy ring lattices with $N=400$ and various noisiness values $\alpha \in [0,1]$. Parameter combinations $(\alpha, T)$ with large ring stability $\Delta$ closely match the regime in which the bifurcation analysis predicts exclusively WFP. This indicates that in regimes which are dominated by WFP, the truncated contagion map consistently recovers the underlying topology of the unit circle, whereas in other parameter regimes this is not possible.

\section{Computational complexity of (truncated) contagion maps }
\label{sec:complexity}

 \begin{figure}[t]{
      \centering 
  \includegraphics[width=0.48\textwidth]{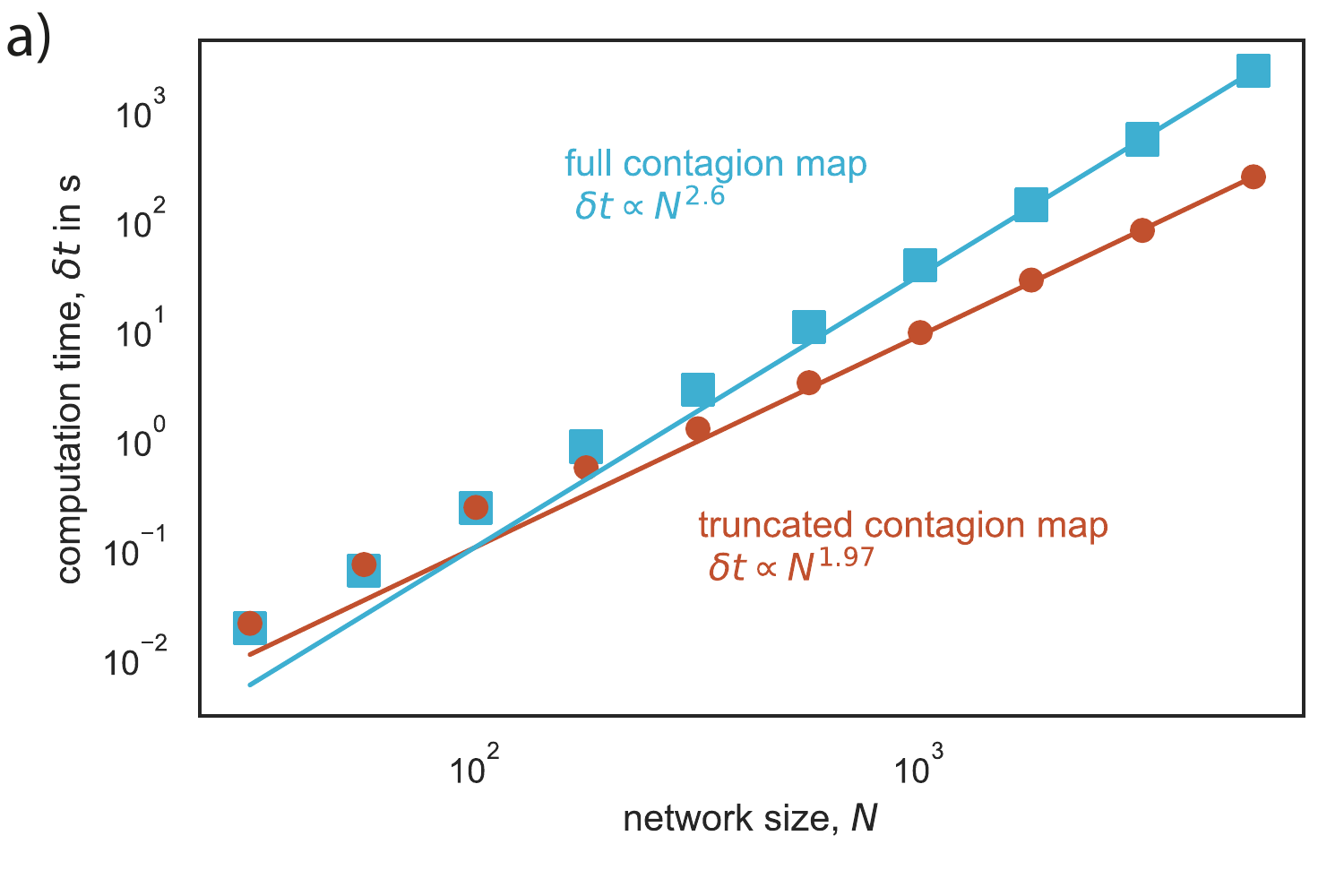}
    \includegraphics[width=0.48\textwidth]{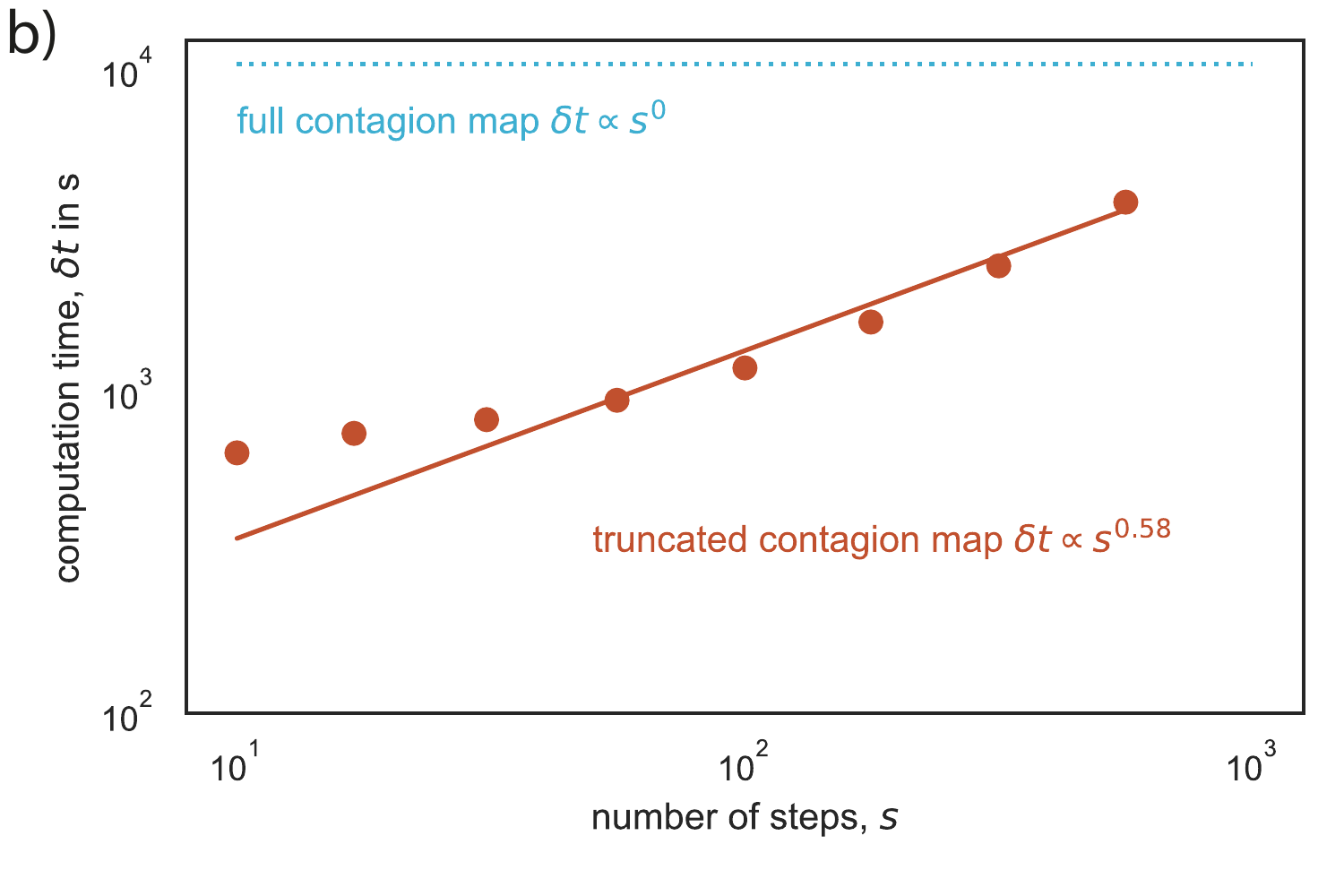}
  }
  \caption{\textbf{Truncated contagion maps have a smaller computational complexity than full contagion maps.} (a) We show the computation time $\delta t$ for contagion maps of noisy ring lattices of various sizes $N \in \{32,10^4\}$, noisiness $\alpha=d^{(NG)} / d^{(G)}= 6/2$, threshold $T=0.3$. The truncated contagion map (red disks; $s=20$ steps considered) is for small network sizes $N<100$ of similar computational cost as the full contagion map (blue squares) but is considerably faster for large networks. We estimate the computational complexity with least-squares fits of a linear equation to the logarithm of the computation time. The computation time $\delta t$ for full and truncated contagion maps scale super-quadratic ($\gamma \approx 2.24$) and sub-quadratic ($\gamma \approx 1.79$) with the network size $N$, respectively. (b) We show the computation time $\delta t$ for contagion maps of noisy ring lattices of size $N =10,000$, noisiness $\alpha=d^{(NG)} / d^{(G)}= 6/2$, and threshold $T=0.3$. We vary the number $s$ of steps and identify an approximately $\gamma \approx 0.58$, with truncated contagion maps approaching the computational time of full contagion maps for $s\to N$.
  }
  \label{fig:complexity}
\end{figure}

\citeauthor{taylor2015topological}\cite{taylor2015topological} established that the typical computational cost of full contagion maps is $\mathcal{O}(N^2d)$, where $N$ is the number of nodes in the input network and $d$ the average degree of its nodes~. As truncated contagion maps omit some activation steps, we anticipate that less computational time is needed than for the full contagion map. 

We empirically verify the reduction of computational cost for contagion maps with thresholds $T=0.3$ on noisy ring lattices with noisiness $\alpha=1/3$ and consider $s=20$ steps in the truncated contagion maps, which is a parameter combination which yields a good recovery of the ring manifold (see Sec.~\ref{sec:topo}). We vary the network sizes $N \in \{32,\dots,10^4 \}$ to identify the computational complexity (see Fig.~\ref{fig:complexity}a). We find that for small network sizes $N\leq 100$, the computation time $\delta t$ of truncated and full contagion maps are indistinguishable because a steady state is reached before a truncation after $s=20$ influences the contagions. Yet, for networks of sizes $N> 100$, the truncated contagion map is quicker to compute than the full contagion map. Furthermore, we identify that the gap between computational costs widens with increasing network size $N$. To obtain estimates of the computational complexity of the form $\delta t = \zeta N^{\gamma}$, we fit its linearised form $\log (\delta t(\log(N))) = \log{\zeta} + {\gamma} \log(N) $ to the logarithm of the data. We find that the computational time of the full contagion map and the truncated contagion map scale superquadratic ($\delta t \propto N^{ 2.24}$) and subquadratic ($\delta t \propto N^{ 1.79}$) with the problem size $N$, respectively, yielding an improved scaling of almost $\sqrt{N}$ for the truncated contagion maps. The reduced computational complexity of the truncated contagion maps can be understood by investigating the time until a steady state is reached. If comparing networks of the same class (i.e., similar characteristics, such as characteristic path length and mean degree), in larger networks the time until a global cascade has been reached tends to be longer than in smaller networks. Thus, a restriction to a fixed number $s$ of steps by the truncated contagion maps, omits a larger number of steps for larger networks. This increasing omission reduces the computational time  for larger networks, leading to a substantial reduction of the computation time for truncated contagion maps.

We further investigate the influence of the number $s$ of steps on the computational time $\delta t$ in Fig.~\ref{fig:complexity}b. Empirically, we find that the computational time grows sublinear with $s$ such that $\delta t \propto s^{0.58}$. As expected, for large number $s$ of steps, the computational time of the truncated contagion maps approaches the computational time of the full contagion maps, highlighting the need for truncation to reduce the computational complexity.

\section{Contagion maps of single-cell transcriptomics data}
\label{sec:singleCell}

\begin{figure*}[t]{
      \centering 
  \includegraphics[width=0.98\textwidth]{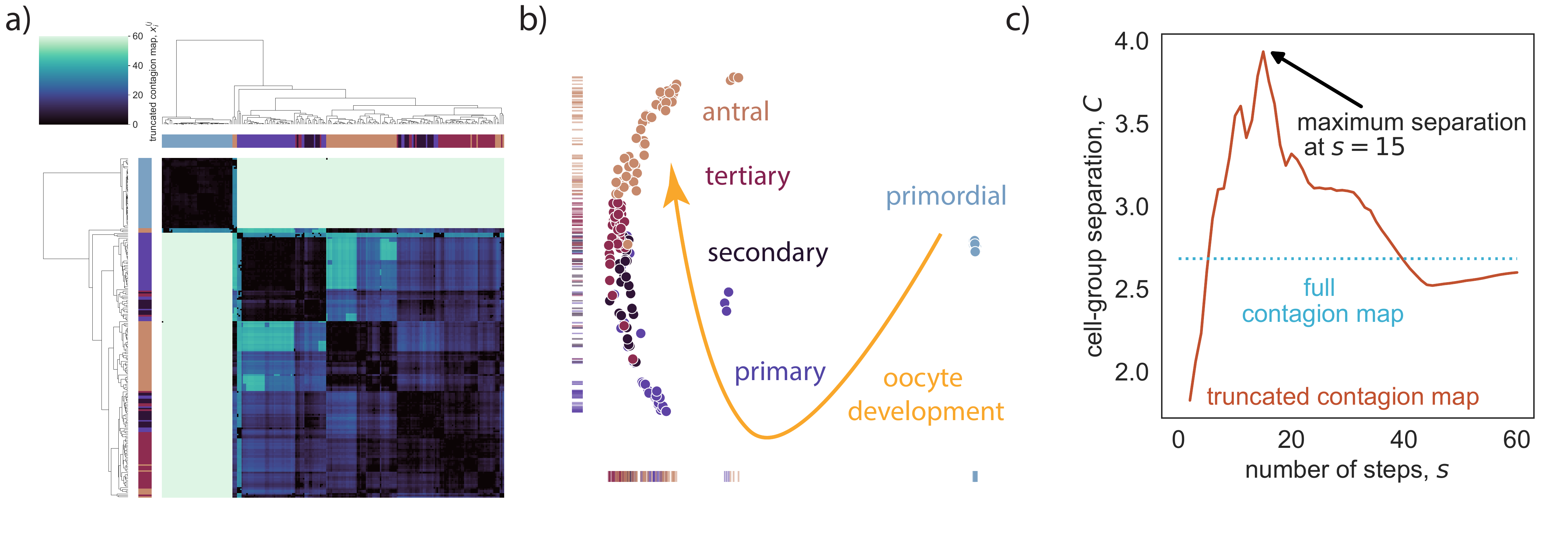}
  }
  \caption{\textbf{Truncated contagion map uncovers differentiation trajectory in single-cell transcriptomics data of mouse oocytes.} (a) We show a truncated contagion map $\mathbf{x}^{(i)}$ with threshold $T=0.3$ and truncation $s=10$ for a cell-similarity network obtained from single-cell transcriptomics data~\cite{gu2019integrative}. The UPGMA method with Euclidean distance was used for the hierarchical clustering. (b) A two-dimensional embedding of the truncated contagion map reveals the oocyte developmental trajectory from primordial oocytes over primary, secondary, and tertiary oocytes to antral oocytes (colour indicates cell identity). (c) The cell-group separation $C$ is highest for a truncation parameter $s=15$ and outperforms a full contagion map (dashed, horizontal line). As expected, for large number $s$ of steps, the performance of the truncated contagion map approaches the performance of the full contagion map.
  }
  \label{fig:singleCell}
\end{figure*}

In this section, we investigate a single-cell transcriptomics dataset of growing mouse oocytes~\cite{gu2019integrative,zhang2021mitochondrial}, in which the oocyte undergo several key stages to generate mature oocytes. Specifically, \citeauthor{gu2019integrative} measure gene expression in five successive stages: primordial follicle, primary oocytes, secondary oocytes, tertiary oocytes, and antral oocytes. The raw gene expression matrix $\mathbf{Count} \in \mathbb{N}^{(c \times g)}$ in which entry $\mathbf{Count}_{i,j}$ indicates the strength of expression of gene $i$ in cell $j$ is available from GEO under GSE114822. For further analysis, we use {\sc scanpy}\cite{wolf2018scanpy} and follow current best practices~\cite{luecken2019current}. For example, we filter out 5937 genes that are detected in less than 3 cells (for details, see online material).  We construct a cell-similarity graph as a $k$-nearest-neighbour graph $G=(V,E)$ in which the $N= 223$ nodes represent the cells $V \in \{1,\dots, c\}$ and each node is connected to its $k=20$ nearest neighbours in gene-expression space.

Then, we construct a truncated contagion map with threshold $T=0.3$ and truncation $s=15$ from the cell-similarity graph (see~Fig.~\ref{fig:singleCell}). The contagion map $\mathbf{x}^{(i)}$ in Fig.~\ref{fig:singleCell}a has a block-diagonal structure in which the the non-growing oocytes (i.e., primordial) are clearly separated from the growing oocyte (primary, secondary, and tertiary) and the fully-grown oocyte (antral). A two-dimensional PCA projection of the contagion map $\mathbf{x}^{(i)}$ in Fig.~\ref{fig:singleCell}b revels the oocyte development trajectory in gene-expression space. The primary oocytes are developing into secondary and tertiary oocytes and finally in antral oocytes. The separation of the developmental stages is not sharp but rather a continuum, yet a pseudo-temporal ordering that matches the cell identities is observable. 

To investigate the performance of the truncated contagion maps in separating the different cell types, we compute the \emph{cell-group separation} $C$, as the ration of median distance between cells of different cell type to the median distance between cells of the same cell type. We investigate this ration $C$ in dependence of the number $s$ of steps in Fig.~\ref{fig:singleCell}c. We find that for most truncation parameters $s$, the cell-group separation $C$ is improved for the truncated contagion maps versus the full contagion map (dashed line). The maximum separation (i.e., best performance) is reached for $s=15$. We note that this optimal choice is close to the heuristic choice of $s^{\ast}=N/10 \approx 22$.

Our study of the single-cell transcriptomics data shows that contagion maps can reveal a manifold structure in a high-dimensional gene expression space by constructing a cell-similarity graph. In particular, we observe that primordial oocytes are clearly separated from developing oocytes and within the developing oocytes, we observe an ordering in gene-expression space.

\section{Conclusion}
\label{sec:conclusion}

In this work, we introduced truncated contagion maps, an extension of contagion maps that has improved computational complexity and is thus computationally less expensive. This is achieved by truncating (i.e., stopping) the contagion processes before they reach a steady state, a technique that is well-established in the context of random walks~\cite{tong2006fast}. Furthermore, we demonstrate by quantifying the persistent homology that, for appropriate parameter choices, truncated contagion maps might improve the recovery of the embedding manifold of noisy geometric networks. We explain this behaviour with the dominance of spreading along the manifold in early steps of a contagion, whereas at later steps the spread along non-geometric edges increases. Lastly, we demonstrated that contagion maps are a manifold-learning technique that can reveal low-dimensional structure in a cell similarity network constructed from single-cell transcriptomics data.


Our findings suggest various avenues for further studies. For example, it is an open question whether these truncation approaches can also be extended to the simplical complex variant of the contagion maps~\cite{kilic2022simplicial}, which is particularly relevant as higher-order interactions are important in cellular interactions~\cite{klamt2009hypergraphs,klimm2021hypergraphs}. Furthermore, we hypothesise that a subsampling approach with selected seeding, as used in the {\tt k-means++} clustering algorithm~\cite{vassilvitskii2006k}, or landmark selection, as used in {\tt L-Isomap}~\cite{de2002global}, might further reduce the computational complexity of contagion maps. Furthermore, the establishment of more sophisticated heuristics for the choice of the truncation parameter $s$ will be fruitful, as its optimal choice is likely to depend on the structure of the analysed network in a complex, non-linear way, which is not reflected in the simple heuristic we proposed in this paper.

\begin{acknowledgments}

FK is supported as an Add-on Fellow for Interdisciplinary Life Science by the Joachim Herz Stiftung and by the Max Planck Society. We thank Martin Vingron for fruitful discussions.
\end{acknowledgments}

\section*{Data Availability Statement}

Data openly available in a public repository that does not issue DOIs. {\sc Python} code to reproduce the results in this paper are available under \url{https://github.com/floklimm/contagionMap}.

\appendix

\bibliography{cmap}

\end{document}